\RequirePackage{ifpdf}
\ifpdf 
\documentclass[pdftex]{sigma}
\else
\documentclass{sigma}
\fi

\begin{document}
\allowdisplaybreaks

\renewcommand{\thefootnote}{$\star$}

\renewcommand{\PaperNumber}{007}

\FirstPageHeading

\ShortArticleName{Quantum Isometry Group for Spectral Triples with Real Structure}

\ArticleName{Quantum Isometry Group for Spectral Triples \\ with Real Structure\footnote{This paper is a
contribution to the Special Issue ``Noncommutative Spaces and Fields''. The
full collection is available at
\href{http://www.emis.de/journals/SIGMA/noncommutative.html}{http://www.emis.de/journals/SIGMA/noncommutative.html}}}

\Author{Debashish GOSWAMI}

\AuthorNameForHeading{D.~Goswami}

\Address{Stat-Math Unit, Indian Statistical Institute, 203, B. T. Road, Kolkata 700108, India}
\Email{\href{mailto:goswamid@isical.ac.in}{goswamid@isical.ac.in}}

\ArticleDates{Received November 06, 2009, in f\/inal form January 17, 2010;  Published online January 20, 2010}

\Abstract{Given a spectral triple of compact type with a real structure in the sense of
 [D\c{a}browski L., \emph{J.~Geom.~Phys.} {\bf 56} (2006), 86--107] (which is a modif\/ication of Connes' original def\/inition to accommodate examples coming from quantum group theory) and refe\-ren\-ces therein,
we prove that there is always  a universal object in the category of compact quantum group acting by orientation preserving isometries (in the sense of [Bhowmick J., Goswami D.,
\emph{J.~Funct.~Anal.} {\bf 257} (2009), 2530--2572]) and also preserving the real structure of the spectral triple. This gives a natural def\/inition of quantum isometry group in the context of real spectral triples without f\/ixing a choice of `volume form' as in [Bhowmick J., Goswami D.,
\emph{J.~Funct.~Anal.} {\bf 257} (2009), 2530--2572].}

\Keywords{quantum isometry groups, spectral triples, real structures}

\Classification{58B32}

\section{Introduction}
Taking motivation from the work of Wang, Banica, Bichon and others (see
\cite{free,wang,ban1,ban2, bichon,univ1} and references therein), we have embarked on a programme to formulate and study various types of `quantum isometry groups' in the setting of (possibly noncommutative) Riemannian geometry. It began with our formulation
 of quantum isometry group based on a `Laplacian' in~\cite{goswami}, and then followed up by a formulation of `quantum group of orientation preserving isometries' in \cite{qorient} (see also \cite{qdisc, q_sphere, jyotish,jyotish_1} for many explicit computations).
The basic idea in all these papers is the following:  f\/irst get an operator theoretic characterisation of an isometric (or orientation preserving and isometric) group action on a Riemannian manifold,  then give an analogous  def\/inition of (compact) quantum group action, and f\/inally try to see whether the category of the compact quantum groups having such action admits a universal object.
However, the transition from group to quantum group action creates a crucial problem, which stems from the fact that unlike the classical group actions  implemented by some unitary representation on a Hilbert space which always preserve the usual trace, a quantum group action may not do so. This problem shows up even in the context of f\/inite dimensional algebras like~$M_n$, and we do not in general get a universal object in the category of quantum groups mentioned before. To get rid of this problem one has to f\/ix a suitable functional (to be interpreted as a choice of `volume form') on the underlying algebra, and then look at the subcategory of the (isometric and orientation preserving) isometric quantum group actions which also preserve this given functional. It has been shown in~\cite{qorient} that this subcategory always has a universal object, which was called there the quantum group of orientation and volume preserving isometries.

The aim of this paper is to provide an alternative to the choice of a volume form. We prove here that if the manifold (possibly noncommutative, i.e.\
given by a spectral triple) has a real structure, then one can get a universal object in the natural subcategory of compact quantum groups whose action, besides being  `orientation-preserving' in the sense of \cite{qorient},  preserves also the real structure in a suitable sense. The idea of the proof is very similar to that of
\cite{qorient}, and we mainly sketch in the present article the arguments which are dif\/ferent from those of \cite{qorient}, but avoid repetition of
  those which are more or
 less the same. The main idea is to  construct a~canonical compact quantum group, which is a free product of countably inf\/initely many copies of the universal
  quantum groups of the form $A_u(Q)$ (notation as in~\cite{wang}), such that any quantum group in the category under consideration can be identif\/ied with a quantum subgroup of this free product.
 In \cite{qorient}, the volume preserving property was used precisely at this step: namely to show that given any eigenvalue $\lambda$ of $D$ there is
a canonical quantum group $A_u(Q_\lambda)$, say, such that the restriction of the action of any quantum group in the above-mentioned category must factor
 through the canonical representation of $A_u(Q_\lambda)$ on the eigenspace corresponding to $\lambda$. The present work relies on the crucial observation
 that a  canonical
 choice (but dif\/ferent from those in \cite{qorient}) of $A_u(Q_\lambda)$ can also be made using the assumption of preservation of the
real structure instead of the volume form
 for an orientation preserving isometric quantum group action.

\section{Preliminaries}
We shall mostly use the notation and terminologies of \cite{qorient}, some of which we brief\/ly recall here again.
We begin by   recalling the def\/inition of compact quantum groups and their actions from   \cite{woro,woro1,vandaelenotes}.  A
compact quantum group (to be abbreviated as CQG from now on)  is given by a pair $({\cal S}, \Delta)$, where ${\cal S}$ is a unital $C^*$-algebra
equipped
 with a unital $\ast$-homomorphism $\Delta : {\cal S} \rightarrow {\cal S} \otimes {\cal S}$ (where $\otimes$ denotes the injective tensor product of $C^*$-algebras)
  satisfying

  $(ai)$ $(\Delta \otimes {\rm id}) \circ \Delta=({\rm id} \otimes \Delta) \circ \Delta$ (co-associativity), and

  $(aii)$ each of the linear spans of $ \Delta({\cal S})({\cal S} \otimes 1)$ and $\Delta({\cal S})(1 \otimes {\cal S})$ is norm-dense in ${\cal S} \otimes {\cal S}$.

      We say that  the compact quantum group $({\cal S},\Delta)$ (co)-acts on a unital $C^*$-algebra ${\cal B}$,
    if there is a  unital $\ast$-homomorphism (called an action) $\alpha : {\cal B} \rightarrow {\cal B} \otimes {\cal S}$ satisfying the following

    $(bi)$ $(\alpha \otimes {\rm id}) \circ \alpha=({\rm id} \otimes \Delta) \circ \alpha$, and

    $(bii)$ the linear span of $\alpha({\cal B})(1 \otimes {\cal S})$ is norm-dense in ${\cal B} \otimes {\cal S}$.

  \begin{definition}
   A unitary (co)representation of a compact quantum group $ ({\cal S},\! \Delta ) $ on a Hilbert~spa\-ce $ {\cal H} $ is a map $ U $ from $ {\cal H} $ to the Hilbert ${\cal S}$-module $ {\cal H} \otimes {\cal S} $  such that the  element $ \widetilde{U} \in {\cal M} ( {\cal K} ( {\cal H} ) \otimes {\cal S} ) $ given by $\widetilde{U}( \xi \otimes b)=U(\xi)(1 \otimes b)$ ($\xi \in {\cal H}, b \in {\cal S})$) is a unitary satisfying  $({\rm  id} \otimes \Delta ) \widetilde{U} = {\widetilde{U}}_{12} {\widetilde{U}}_{13},$ where for an operator $X \in {\cal B}({\cal H}_1 \otimes {\cal H}_2)$ we have denoted by $X_{12}$ and $X_{13}$ the operators $X \otimes I_{{\cal H}_2} \in {\cal B}({\cal H}_1 \otimes {\cal H}_2 \otimes {\cal H}_2)$, and $\Sigma_{23} X_{12} \Sigma_{23}$ respectively ($\Sigma_{23}$ being the unitary on ${\cal H}_1 \otimes {\cal H}_2 \otimes {\cal H}_2$ which f\/lips the two copies of ${\cal H}_2$).

Given a unitary representation $U$ we shall denote by $\alpha_U$ the $\ast$-homomorphism $\alpha_U(X)=\widetilde{U}(X \otimes 1){\widetilde{U}}^*$ for $X \in {\cal B}({\cal H})$. For a  not necessarily bounded, densely def\/ined (in the weak operator topology)  linear functional $\tau$ on ${\cal B}({\cal H})$,  we say that $\alpha_U$ preserves $\tau$ if $\alpha_U$ maps a suitable (weakly) dense $\ast$-subalgebra   (say ${\cal D}$) in the domain of $\tau$ into ${\cal D} \otimes_{\rm alg} {\cal S}$ and $( \tau \otimes {\rm id}) (\alpha_U(a))=\tau(a)1_{\cal S}$  for all $a \in {\cal D}$. When $\tau$ is bounded and normal, this is equivalent to $(\tau \otimes {\rm id}) (\alpha_U(a))=\tau(a) 1_{\cal S}$ for all $a \in {\cal B}({\cal H})$.

We say that a (possibly unbounded) operator $T$ on ${\cal H}$ commutes with $U$ if $T \otimes I$ (with the natural domain) commutes with $\widetilde{U}$. Sometimes such an operator will be called $U$-equivariant.
\end{definition}

 Let us now recall the concept of universal quantum groups as in
\cite{univ1,free}
and references therein. We shall use most of the terminologies of
\cite{free}, e.g.\ Woronowicz $C^*$-subalgebra, Woronowicz
$C^*$-ideal etc, however with the exception that we shall call the
Woronowicz $C^*$-algebras just compact quantum groups, and not use
the term compact quantum groups for the dual objects as done in~\cite{free}.
For an $n \times n$ positive invertible matrix $Q=((Q_{ij}))$, let $A_u(Q)$ be the compact quantum group def\/ined and  studied in \cite{wang,univ1}, which is the universal $C^{*}$-algebra generated by $ \{ u^{Q}_{kj}, \ k,j=1,\dots ,n \}$ such that $u:=((u_{kj} \equiv  u^{Q}_{kj} ))$ satisf\/ies \begin{equation} \label{wangalg} u u^*=I_n =u^{*}u, \qquad u^{\prime} Q  \overline{u} Q^{-1}=I_n=Q{\overline{u}} Q^{-1} u^{\prime}.\end{equation} Here
$u^{\prime} =(( u_{ji} ))$ and $\overline{u}=(( u_{ij}^{*} ))$, and also note that we have made the identif\/ication of an $n \times n$ matrix $B$ with its trivial ampliation $B \otimes 1$ in $M_n({\mathbb C}) \otimes {\cal A}$ for any $C^*$-algebra ${\cal A}$. The coproduct, say $\tilde{\Delta}$, is given by,  $\tilde{\Delta}(u_{ij})=\sum_k u_{ik} \otimes u_{kj}.$ It may be noted that $A_u(Q)$ is the universal object in the category of compact quantum groups generated by the coef\/f\/icients
of a unitary representation~$v$ on~${\mathbb C}^n$ such that the
adjoint action ${\rm Ad}_v$ on $M_n({\mathbb C})$  preserves  the functional $M_n \ni x \mapsto {\rm Tr({Q}^{\prime} x)}$ (see~\cite{wangergodic}), where we
refer the reader to \cite{univ1}  for a detailed  discussion on the structure and classif\/ication of
such quantum groups.

Given a  $C^*$-algebra ${\cal S}$ we shall denote by $\tilde{J}_{\cal S}$ the antilinear map $a \mapsto a^*$. For any faithful state
(which exists whenever ${\cal S}$ is separable) this map can be viewed as a closable unbounded antilinear map on the GNS space of the state,
and the corresponding closed extension will be denoted by the same notation.

We now give a def\/inition of the real structure along the lines of~\cite{dab_real} and~\cite{landi_real}, which is a suitable modif\/ication of Connes' original  def\/inition (see~\cite{connes}) to accommodate the examples coming from quantum groups and quantum homogeneous spaces.

\begin{definition}
An odd spectral triple  with a real structure is given by a spectral triple $({\cal A}^\infty{,} {\cal H}{,} D)$ along with a (possibly unbounded, invertible)
closed anti-linear operator $\tilde{J}$ on ${\cal H}$ such that  ${\cal D}:={\rm Dom}(D) \subseteq {\rm Dom}(\tilde{J}) $,
$\tilde{J} {\cal D} \subseteq {\cal D}$,  $\tilde{J}$ commutes with $D$ on ${\cal D}$, and the antilinear isometry $J$ obtained from the polar decomposition of $\tilde{J}$
 satisf\/ies the usual conditions for a real structure in the sense of~\cite{landi_real}, for a suitable sign-convention given by $(\epsilon, \epsilon^\prime)
  \in \{ \pm 1 \} \times \{ \pm 1 \}$ as described in \cite[page~30]{varilly}, i.e.\  $J^2=\epsilon I$, $JD=\epsilon^\prime DJ$, and for all  $x,y \in {\cal A}^\infty$, the commutators $[x, JyJ^{-1}]$ and $[JxJ^{-1},[D,y]]$ are compact operators.

If the spectral triple is even,  a real structure with the sign-convention given by a triplet $(\epsilon, \epsilon^\prime, \epsilon^{\prime \prime})$   as in~\cite[page~30]{varilly} is similar to a real structure in the odd case (with the sign-convention $(\epsilon, \epsilon^\prime)$), but with the  additional requirement that  $J \gamma = \epsilon^{\prime \prime} \gamma J$.
\end{definition}

We now recall from \cite{qorient} the def\/inition of quantum family of orientation preserving isometries and then appropriately adapt it to the framework of real structure.

\begin{definition}
         \label{def_q_fam} {\rm A quantum family of orientation  preserving  isometries for the spectral triple $({{\cal A}^\infty}, {\cal H}, D)$ is given by a pair $({\cal S}, U)$ where ${\cal S}$ is a separable unital $C^*$-algebra and  $U$ is an ${\mathbb C}$-linear map from ${\cal H}$ to the Hilbert module ${\cal H} \otimes {\cal S}$ such that the ${\cal S}$-linear map $\widetilde{U}$ given by $\widetilde{U}( \xi \otimes b)=U(\xi) (1 \otimes b)$ $(\xi \in {\cal H}$, $b \in {\cal S}$) extends to a unitary element of  $ {\cal M}({\cal K}({\cal H}) \otimes {\cal S})$ satisfying the following

$(i)$ $\tilde{U}$ commutes with $D \otimes I$, and

$(ii)$ $({\rm id} \otimes \phi) \circ \alpha_U(a) \in ({{\cal A}^\infty})^{\prime \prime}$ $\forall a \in {\cal A}^\infty$ for every state $\phi$ on ${\cal S}$, where $\alpha_U(x):=\widetilde{U}( x \otimes 1) {\widetilde{U}}^* $ for $x \in {\cal B}({\cal H})$.

In case the $C^*$-algebra ${\cal S}$ has a coproduct $\Delta$ such that $({\cal S},\Delta)$ is a compact quantum group and~$U$ is a unitary representation  of~$({\cal S}, \Delta)$ on~${\cal H}$, we say that $({\cal S}, \Delta)$ acts  by orientation preserving isometries  on the spectral triple.}
\end{definition}

Given a quantum family of orientation preserving isometries $({\cal S}, U)$ as above, note that, since~$D$ has f\/inite dimensional eigenspaces which are preserved by~$U$, we have $U{\cal D}_0 \subseteq {\cal D}_0 \otimes_{\rm alg} {\cal S}$, where~${\cal D}_0$ denotes the linear span of eigenvectors of~$D$.

\begin{definition}
Suppose that the (odd) spectral triple $({\cal A}^\infty, {\cal H}, D)$ is equipped with a real structure given by $\tilde{J}$. We say that a quantum family of orientation preserving isometries $({\cal S}, U)$ also preserves the real structure if the following holds on ${\cal D}_0$:
\begin{equation} \label{equiv_real}(\tilde{J} \otimes \tilde{J}_{\cal S}) \circ U=U \circ \tilde{J}.\end{equation}
In case the $C^*$-algebra ${\cal S}$ has a coproduct $\Delta$ such that $({\cal S},\Delta)$ is a compact quantum group and~$U$ is a unitary representation  of $({\cal S}, \Delta)$ on~${\cal H}$, we say that $({\cal S}, \Delta)$ acts  by orientation and real structure  preserving isometries  on the spectral triple.

Similar def\/initions can be given in the even case, with the additional requirement being that~$U$ commutes with $\gamma$.
\end{definition}

Given a compact quantum group ${\cal Q}$ acting on ${\cal A}$, such that the action is implemented by a unitary representation $U$ of the quantum group on ${\cal H}$,
  it is easy to see that the notion of equivariance of the spectral triple with the real structure as proposed in  \cite{dab_real} is equivalent to saying that $({\cal Q}, U)$ is a quantum group acting by orientation and real structure preserving isometries in our sense. We refer the reader to~\cite{dab_real} for related discussions and examples of such equivariant real spectral triples.

As in~\cite{qorient}, we consider the category ${\bf Q} \equiv {\bf Q}(D)$ with the object-class consisting of all quantum families of orientation and real structure  preserving isometries $({\cal S}, U)$ of the given spectral triple, and the set of morphisms ${\rm Mor}(({\cal S},U),({\cal S}^\prime,U^\prime))$ being the set of unital $\ast$-homomorphisms $\Phi : {\cal S} \rightarrow {\cal S}^\prime$ satisfying $({\rm id} \otimes \Phi) (U)=U^\prime$. We also consider another category ${\bf Q}^\prime \equiv {\bf Q}^\prime(D)$ whose objects are triplets $({\cal S}, \Delta, U)$,  where $({\cal S},\Delta)$ is a compact quantum group acting by orientation and real structure preserving isometries on the given spectral triple, with $U$ being the corresponding unitary representation. The morphisms  are the homomorphisms of compact quantum groups which are also morphisms of the underlying quantum families of orientation preserving isometries. The forgetful functor $F: {\bf Q}^\prime \rightarrow {\bf Q}$ is clearly faithful, and we can view $F({\bf Q}^\prime)$ as a~subcategory of ${\bf Q}$.
Our aim is to show that the above categories admit universal object, which we prove in the next section.

\section{Main results and examples}

Let us  f\/ix a spectral triple $({\cal A}^\infty, {\cal H}, D)$ which is of compact type along with a real structure given by $\tilde{J}$. We shall work with an odd spectral triple, but remark that all the arguments will go through almost verbatim, with some obvious and minor changes at places, in the even case.  The sign-convention of the real structure is not explicitly mentioned, since it is not going to be needed anywhere, and we remark that our arguments are valid for any possible choice of the signs. The $C^*$-algebra generated by ${\cal A}^\infty$ in ${\cal B}({\cal H})$ will be denoted by ${\cal A}$. Let $\lambda_0=0, \lambda_1, \lambda_2, \ldots$ be the eigenvalues of $D$ with $V_i$ denoting the ($d_i$-dimensional, $d_i<\infty$) eigenspace for~$\lambda_i$. Let $\{ e_{ij}, j=1,\dots , d_i \}$ be an orthonormal basis of $V_i$. Clearly, $\{ \tilde{J}(e_{ij}), \ i \geq 0, \ 1 \leq j \leq d_i \}$ is a linearly independent (but not necessarily orthogonal) set, and let $T_i$ denote the positive nonsingular matrix $\big( \langle  \tilde{J}(e_{ij}), \tilde{J}(e_{ik}) \rangle \big)_{j,k=1}^{d_i}.$
 Let us denote the CQG $A_u(T_i)$ by ${\cal U}_i$, with its canonical unitary representation $\beta_i$ on $V_i \cong {\mathbb C}^{d_i}$, given by $\beta_i(e_{ij})=\sum_k e_{ik} \otimes u^{T_i}_{kj}$. Let ${\cal U}$ be the free product of ${\cal U}_i$, $i=1,2,\dots $ and $\beta=\ast_i \beta_i$ be the corresponding free product representation of ${\cal U}$ on ${\cal H}$. We shall also consider the  corresponding  unitary element $\tilde{\beta}$ in ${\cal M}({\cal K}({\cal H}) \otimes {\cal U})$.

 \begin{lemma}
\label{lem2} Consider the real spectral triple $({\cal A}^\infty,{\cal H},D, \tilde{J})$ as before   and let $({\cal S},U)$ be a quantum family of  orientation and real structure  preserving  isometries of the given spectral triple.  Moreover, assume that the map  $U$ is faithful in the sense that there is no
proper  $C^*$-subalgebra~${\cal S}_1$ of~${\cal S}$ such that
$\widetilde{U} \in {\cal M}({\cal K}({\cal H} ) \otimes {\cal S}_1)$.
Then  we
can find a $\ast$-isomorphism  $\phi : {\cal U}/
I \rightarrow {\cal S}$ between~${\cal S}$ and a~quotient of~${\cal U}$ by a
 $C^*$-ideal
 $I$ of ${\cal U}$, such that $ U= ({\rm id}\otimes \phi) \circ ({\rm id} \otimes \Pi_I) \circ
 \beta$, where~$\Pi_I$ denotes the quotient map from~${\cal U}$ to
 ${\cal U}/I$.

 If, furthermore, there is a compact quantum group structure on ${\cal S}$ given by a coproduct  $\Delta$ such that $({\cal S},\Delta, U)$ is an object in ${\bf Q}^\prime(D)$, the ideal $I$ is a Woronowicz $C^*$-ideal and the $\ast$-isomorphism $\phi : {\cal U}/ I \rightarrow {\cal S}$ is a morphism of compact quantum groups.
\end{lemma}

\begin{proof}
We follow the line of arguments of a similar result in \cite{qorient}, though with suitable mo\-di\-f\/i\-ca\-tions.
It is clear that $U$  maps $V_i $ into $V_i \otimes
{\cal S}$ for each $i$. Let $v^{(i)}_{kj}$ ($j,k=1,\dots ,d_i$) be the
elements of ${\cal S}$ such that $U(e_{ij})=\sum_k e_{ik} \otimes
v^{(i)}_{kj}$.  Note that $v_i:=((v^{(i)}_{kj} ))$ is a unitary in
$M_{d_i}({\mathbb C}) \otimes {\cal S}$. Moreover, the $\ast$-subalgebra
generated by all $ \{ v^{(i)}_{kj}, \ i \geq 0, \ j,k=1 ,\dots ,d_i\}$ must be
dense in ${\cal S}$ by the assumption of faithfulness.

Now, we shall make use of~(\ref{equiv_real}). Fix any $i$ and let $\Lambda_i=(( \tau_{lm} ))$ be the matrix such that $\tilde{J}(e_{ij})=\sum_l \tau_{lj} e_{il}$. By assumption, $\Lambda_i$ is invertible, and it is clear that $\Lambda_i^* \Lambda_i=T_i$.
Expanding both sides of
$U(\tilde{J}e_{ij})=\sum_k \tilde{J}e_{ik} \otimes (v^{(i)}_{kj})^*$ we get \begin{equation} \label{123}\sum_m e_{im} \otimes \left( \sum_l \tau_{lj} v^{(i)}_{ml} \right)=
 \sum_m e_{im} \otimes \left( \sum_l \tau_{ml} (v^{(i)}_{lj})^* \right).\end{equation}  By comparing coef\/f\/icients of $e_{im}$ in both sides of  (\ref{123}), we get
 $\sum_l \tau_{lj} v^{(i)}_{ml}=\sum_l \tau_{ml} (v^{(i)}_{lj})^*,$ that is, $v_i \Lambda_i=\Lambda_i \overline{v_i}$. It follows that
  $\overline{v_i}=\Lambda_i^{-1} v_i \Lambda_i, $ hence $\overline{v_i}$ is invertible, since $v_i$ is so.
Moreover, taking the ${\cal S}$-valued inner product $\langle \cdot, \cdot \rangle_{\cal S}$ on both sides of $U(\tilde{J}e_{ij})=\sum_k \tilde{J}e_{ik} \otimes (v^{(i)}_{kj})^*$   we obtain $T_i=v_i^\prime T_i \overline{v_i}$. Thus, $T_i^{-1}v_i^\prime T_i$ must be the (both-sided) inverse of $\overline{v_i}$, from which we see that the relations (\ref{wangalg}) are satisf\/ied with $u$ replaced by $v_i$.

We get,  by the universality of ${\cal U}_i$,  a
$\ast$-homomorphism from ${\cal U}_i$ to ${\cal S}$ sending $u^{(i)}_{kj} \equiv u_{kj}^{T_i} $
to $v^{(i)}_{kj}$, and by def\/inition of the free product, this
induces a $\ast$-homomorphism, say $\Pi$, from ${\cal U}$ onto ${\cal S}$,
 so that ${\cal U}/I
\cong {\cal S}$, where  $I:={\rm Ker}(\Pi)$.

In case ${\cal S}$ has a coproduct $\Delta$ making it into a compact quantum group and $U$  is a quantum group representation, it is easy to see that  the subalgebra of ${\cal S}$ generated by $\{ v^{(i)}_{kj},\ i \geq 0, \ j,k=1,\dots ,d_i \}$ is a Hopf
algebra, with $\Delta(v^{(i)}_{kj})=\sum_l v^{(i)}_{kl} \otimes v^{(i)}_{lj}$. From this, it follows that  $\Pi$ is Hopf-algebra morphism, hence $I$ is a Woronowicz $C^*$-ideal.
\end{proof}
\begin{remark}
From the proof of the above result, it can be seen that the assumption of preserving the real structure implies that the `$R$-twisted volume form' is preserved, where $R$ is given by $R|_{V_i}=T_i^\prime$. This connects the approach of the present article to that of \cite{qorient}, and in some sense gives an explanation of how the proof of the above lemma works.
\end{remark}

The rest of the arguments in \cite{qorient} goes through more or less verbatim  and we have the following analogue of the main result of \cite{qorient}:

\begin{theorem}
For any  real $($odd or even$)$ spectral triple $({\cal A}^\infty, {\cal H}, D, \tilde{J})$, the category ${\bf Q}$ of quantum families of orientation and real structure preserving isometries has a universal $($initial$)$ object, say  $(\widetilde{{\cal G}}, U_0)$. Moreover, $\widetilde{{\cal G}}$ has a coproduct $\Delta_0$ such that $(\widetilde{{\cal G}},\Delta_0)$ is a compact quantum group and $(\widetilde{{\cal G}},\Delta_0,U_0)$ is a universal object in the category ${\bf Q}^\prime$.  The representation  $U_0$ is faithful.
\end{theorem}

\begin{definition}
Let ${\cal G}$ denote the Woronowicz subalgebra of $\widetilde{{\cal G}}$ generated by elements of the form $\langle \xi \otimes 1, {\rm ad}_{U_0}(a)(\eta \otimes 1) \rangle $, where $\xi, \eta \in {\cal H}$, $a \in {\cal A}^\infty$, and where $\langle  \cdot, \cdot\rangle$ denotes the $\widetilde{{\cal G}}$-valued inner product of the Hilbert module ${\cal H} \otimes \widetilde{{\cal G}}$. We shall call ${\cal G}$  the quantum  group of orientation and real structure preserving isometries of the given  spectral triple, and denote it by $ {QISO}^+ ( {\cal A}^{\infty}, {\cal H}, D, \tilde{J} ) $ or even simply as ${QISO}^{+}_{{\rm real}}(D)$. The quantum group $\widetilde{{\cal G}}$ is denoted by $\widetilde{QISO}^{+}_{{\rm real}}(D)$.
\end{definition}

\begin{remark}\label{11}
It is clear from the def\/inition that $QISO^+_{{\rm real}}\!(D)$ is a quantum subgroup of $QISO^+\!(D)\!$ whenever the later exists, since the former is the universal object in a subcategory of the category for which the latter is universal (if exists).
\end{remark}

We conclude the article with two examples.

\begin{example}
The standard spectral triple on the noncommutative two torus ${\cal A}_\theta$ (including the commutative case, i.e. $\theta=0$) has a canonical
real structure. The Hilbert space ${\cal H}$ is in this case $L^2({\cal A}_\theta, \tau) \otimes {\mathbb C}^2$ (where $\tau$ is the canonical faithful trace on ${\cal A}_\theta$)
 and $D=\left( \begin{array}{cc} 0 & d_1+id_2\\ d_1-id_2 & 0 \end{array} \right), $
$J=\left( \begin{array}{cc} 1 & 0\\ 0 & -1 \end{array} \right),$ where $d_1$, $d_1$ denote the canonical derivations (see~\cite{connes}).
 This spectral  triple (without taking into account the real structure) has been considered in Subsection~4.3  of~\cite{qorient}, where we have proved that the quantum group orientation preserving isometries exists and coincides with the classical group of such isometries, i.e.~$C({\mathbb T}^2)$. Since it can easily be seen that  this $C({\mathbb T}^2)$ action also preserves the real structure, it follows from Remark~\ref{11} that $QISO^+_{{\rm real}}(D)$ must be $C({\mathbb T}^2)$. \end{example}

\begin{example}
This is an example involving a quantum group action with nontrivial modularity. Consider the spectral triple on the Podles sphere $S^2_{\mu c}$
constructed in~\cite{Dabrowski_et_al}. Note that in~\cite{Dabrowski_et_al}, a real structure has also been constructed and the spectral triple
as well as the real structure are shown to be equivariant in the sense of \cite{dab_real} with respect to the canonical action of $SO_\mu(3)$.
Thus, $QISO^+_{{\rm real}}(D)$ for this spectral triple has $SO_\mu(3)$ as a quantum subgroup, and then it follows from Theorem~3.39  of~\cite{q_sphere} that $QISO^+_{{\rm real}}(D)$ must coincide with $SO_\mu(3)$.
\end{example}

\subsection*{Acknowledgements}
The author acknowledges the support from Indian National Science Academy for the project `Noncommutative Geometry and Quantum Groups' and UKIERI, British Council.

\pdfbookmark[1]{References}{ref}
\LastPageEnding


\begin{thebibliography}{99}

\footnotesize\itemsep=0pt

\bibitem{ban1}
Banica T.,
Quantum automorphism groups of small metric spaces,
\href{http://dx.doi.org/10.2140/pjm.2005.219.27}{\emph{Pacific J. Math.}} {\bf 219} (2005), 27--51,
\href{http://arxiv.org/abs/math.QA/0304025}{math.QA/0304025}.


\bibitem{ban2}
Banica T.,
Quantum automorphism groups of homogeneous graphs,
\href{http://dx.doi.org/10.1016/j.jfa.2004.11.002}{\emph{J.~Funct. Anal.}} {\bf 224} (2005), 243--280,
\href{http://arxiv.org/abs/math.QA/0311402}{math.QA/0311402}.

\bibitem{bichon}
Bichon J.,
Quantum automorphism groups of f\/inite graphs,
\href{http://dx.doi.org/10.1090/S0002-9939-02-06798-9}{\emph{Proc.\ Amer.\ Math.\ Soc.}} {\bf 131} (2003), 665--673,
\href{http://arxiv.org/abs/math.QA/9902029}{math.QA/9902029}.


\bibitem{jyotish_1}
Bhowmick J.,
Quantum isometry group of the $n$-tori,
\href{http://dx.doi.org/10.1090/S0002-9939-09-09908-0}{\emph{Proc. Amer. Math. Soc.}} {\bf 137} (2009), 3155--3161,
\href{http://arxiv.org/abs/0803.4434}{arXiv:0803.4434}.

\bibitem{qorient}
Bhowmick J., Goswami D.,
Quantum group of orientation-preserving Riemannian isometries,
\href{http://dx.doi.org/10.1016/j.jfa.2009.07.006}{\emph{J. Funct. Anal.}} {\bf 257} (2009), 2530--2572,
\href{http://arxiv.org/abs/0806.3687}{arXiv:0806.3687}.

\bibitem{qdisc}
Bhowmick J., Goswami D., Skalski A.,
Quantum isometry groups of 0-dimensional manifolds,
\emph{Trans. Amer. Math. Soc.}, to appear,
\href{http://arxiv.org/abs/0807.4288}{arXiv:0807.4288}.

\bibitem{jyotish}
Bhowmick J., Goswami D.,
Quantum isometry groups: examples and computations,
\href{http://dx.doi.org/10.1007/s00220-008-0611-5}{\emph{Comm. Math. Phys.}} {\bf 285} (2009), 421--444,
\href{http://arxiv.org/abs/0707.2648}{arXiv:0707.2648}.

\bibitem{q_sphere}
Bhowmick J., Goswami D.,
Quantum isometry groups of the Podles sphere,
\href{http://arxiv.org/abs/0810.0658}{arXiv:0810.0658}.


\bibitem {connes}
Connes A.,
Noncommutative geometry, Academic Press, Inc., San Diego, CA, 1994.

\bibitem{dab_real}
D\c{a}browski L.,
Geometry of quantum spheres,
\href{http://dx.doi.org/10.1016/j.geomphys.2005.04.003}{\emph{J.~Geom. Phys.}} {\bf 56} (2006), 86--107,
\href{http://arxiv.org/abs/math.QA/0501240}{math.QA/0501240}.

\bibitem{landi_real}
D\c{a}browski L., Landi G., Paschke M.,  Sitarz A.,
The spectral geometry of the equatorial Podle\'s sphere,
\href{http://dx.doi.org/10.1016/j.crma.2005.04.003}{\emph{C.~R.~Math. Acad. Sci. Paris}} {\bf 340} (2005),  819--822,
\href{http://arxiv.org/abs/math.QA/0408034}{math.QA/0408034}.

\bibitem{Dabrowski_et_al}
D\c{a}browski L., D'Andrea F., Landi G., Wagner E.,
Dirac operators on all Podle\'s quantum spheres,
\emph{J.~Noncommut. Geom.} {\bf 1} (2007), 213--239,
\href{http://arxiv.org/abs/math.QA/0606480}{math.QA/0606480}.

\bibitem{goswami}
Goswami D.,
Quantum group of isometries in classical and noncommutative geometry,
\href{http://dx.doi.org/10.1007/s00220-008-0461-1}{\emph{Comm. Math. Phys.}} {\bf 285} (2009),  141--160,
\href{http://arxiv.org/abs/0704.0041}{arXiv:0704.0041}.

\bibitem{vandaelenotes}
Maes A., Van Daele A.,
Notes on compact quantum groups,
\emph{Nieuw Arch.\ Wisk. (4)} {\bf 16} (1998), 73--112,
\href{http://arxiv.org/abs/math.FA/9803122}{math.FA/9803122}.


\bibitem{varilly}
Varilly J.C.,
An introduction to noncommutative geometry, {\it EMS Series of Lectures in Mathematics}, European Mathematical Society (EMS), Z\"urich, 2006.


\bibitem{free}
Wang S.,
Free products of compact quantum groups,
\href{http://dx.doi.org/10.1007/BF02101540}{\emph{Comm.\ Math.\ Phys.}} {\bf 167} (1995),  671--692.

\bibitem{wang}
Wang S.,
Quantum symmetry groups of f\/inite spaces,
\href{http://dx.doi.org/10.1007/s002200050385}{\emph{Comm.\ Math.\ Phys.}} {\bf 195} (1998), 195--211,
\href{http://arxiv.org/abs/math.OA/9807091}{math.OA/9807091}.

\bibitem{univ1}
Wang S.,
Structure and isomorphism classif\/ication of compact quantum groups $A_u(Q)$ and $B_u(Q)$,
\emph{J.~Ope\-rator Theory} {\bf 48} (2002), 573--583,
\href{http://arxiv.org/abs/math.OA/9807095}{math.OA/9807095}.



\bibitem{wangergodic}
 Wang S.,
 Ergodic actions of universal quantum groups on operator algebras,
\href{http://dx.doi.org/10.1007/s002200050622}{\emph{Comm. Math. Phys.}} {\bf 203} (1999), 481--498,
\href{http://arxiv.org/abs/math.OA/9807093}{math.OA/9807093}.

\bibitem{woro1}
Woronowicz S.L.,
Compact matrix pseudogroups,
\href{http://dx.doi.org/10.1007/BF01219077}{\emph{Comm.\ Math.\ Phys.}} {\bf 111} (1987), 613--665.

\bibitem{woro}
Woronowicz S.L.,
Compact quantum groups, in  Sym\'etries Quantiques
 (Les Houches, 1995), Editors A.~Connes et al., North-Holland, Amsterdam,  1998, 845--884.

\end{thebibliography}
\end{document}